\documentclass[11pt]{article}
\usepackage{amsmath}
\usepackage{amsmath,amssymb,amsfonts,amsthm,fancyhdr}
\usepackage{epsfig,graphicx,picins,picinpar,subfigure}
\usepackage{pstricks}
\usepackage{fancyvrb}
\usepackage[numbers,sort&compress]{natbib}
\begin{document}

\title{\textbf{A Refinement of the Function $g(x)$ on Grimm's Conjecture}}
\author{\emph{Shaohua Zhang}$^{1,2}$}
\date{{\small 1 School of Mathematics, Shandong University,
Jinan,  Shandong, 250100, China \\
2 The key lab of cryptography technology and information security,
Ministry of Education, Shandong University, Jinan, Shandong, 250100,
China\\E-mail address: shaohuazhang@mail.sdu.edu.cn}}
 \maketitle

\begin{abstract}
In this paper, we refine the function $g(x)$ on Grimm's conjecture
and improve a result of Erd\"{o}s and Selfridge without using Hall's
theorem.

\vspace{3mm}\textbf{Keywords:} consecutive composite numbers,
Grimm's Conjecture, Cram\'{e}r's conjecture, binomial coefficient

\vspace{3mm}\textbf{2000 MR  Subject Classification:}\quad 11A41,
11A99, 11B65
\end{abstract}


\section{Some basic notations}
\setcounter{section}{1}\setcounter{equation}{0}

Let $\pi (x)$ be the prime counting function which represents the
number of primes not exceeding the real number $x$. We write $f(x)=O
(g(x))$, or equivalently $f(x)\ll g(x)$ when there is a constant $C$
such that $|f(x)|\leq C g(x)$ for all values of $x$ under
consideration. We write $f(x)=o (g(x))$ when $\lim_{x\rightarrow
\infty} f(x)/g(x)=0$. We let $[x]$ denote the largest integer not
exceeding the real number $x$.

Let $\binom{m+n}{n}$ be the binomial coefficient. For a prime $p$
and a positive integer $n>1$,  we define $v_p(n)$ to be the largest
exponent of $p$ that divides $n$. In this paper, we always denote by
$p$ a prime number. Denote the set of all prime numbers by $P$ and
denote  the set of all composite numbers by $C$. Let $H_n=\{x:x\in
C,\forall p|x, p^{v_p(x)}\leq n \}$.

\section{Introduction}
In 1969, Grimm [2] made an important conjecture that if
$m+1,...,m+n$ are  consecutive composite numbers, then there exist
$n$ distinct prime numbers  $p_1,...,p_n$ such that $m+j$  is
divisible by $p_j$, $1\leq j\leq n$. This implies that the product
of any $n$ consecutive composite numbers must have at least $n$
distinct prime factors.

Grimm proved the conjecture for two special cases: (i) For all $n$
in the sequence of consecutive composites, $\{n!+i\},i=2,...,n$.
(ii) For all $m$, when $m>n^{n-1}$. This was improved to $m>n^{\pi
(n)}$ by Erd\"{o}s and Selfridge [1, Theorem 3] who used Hall's
theorem. Moreover, in [1, Theorem 1], Erd\"{o}s and Selfridge
obtained the following typical theorem: \emph{Let $v(m,n)$  be the
number of distinct prime factors of the product
$\prod_{i=1}^{i=n}(m+i)$, and let $f(m)$ be the largest $n$ for
which $v(m,n)\geq n$, then $f(m)<c(\frac{m}{\log m})^\frac{1}{2}$
for all $m$, where $c$ is a positive absolute constant.}

Thus, Grimm's Conjecture implies $p_{r+1}-p_r\ll (p_r/\log
p_r)^\frac{1}{2}$ which is out of bounds for even the Riemann
Hypothesis [14, Introduction] which implies that $p_{r+1}-p_r\ll
p_r^\frac{1}{2}\log p_r$, where $p_{r+1}$ and $p_r$ are consecutive
primes. It implies particularly that there are primes between  $n^2$
and $(n+1)^2$ for all sufficiently large $n$, a conjecture which is
still open. For the details of its proof, see also [4, Appendix 1].
Furthermore, by Theorem 1 in [1], it is not difficult to prove that
Grimm's Conjecture implies that there are primes between   $x^2$ and
for $x^2+x$  all sufficiently large $x$ and there are primes between
$y^2-y$ and $y^2$ for all sufficiently large $y$. Let $x=n$ and
$y=n+1$, then, Grimm's Conjecture implies that for all sufficiently
large $n$, there are two primes between  $n^2$ and $(n+1)^2$  which
implies that there are four primes between  $n^3$ and $(n+1)^3$  for
all sufficiently large $n$ [5]. These surprising consequences
motivate the study of the function  $g(m)$ which is the largest
integer $n$ such that Grimm's Conjecture holds for the interval
$[m+1,m+n]$. Paulo Ribenboim [4, Appendix 1] pointed out that
$g(m)<2m$ since the interval will contain clearly two powers of 2.
Theorem 1 in [1] said that $g(m)=O(\sqrt{m/\log m})$. By using a
result of Ramachandra [6], Erd\"{o}s and Pomerance [14] pointed out
that $g(m)<m^{\frac{1}{2}-c}$  for some fixed $c>0$ and all large
$m$. Thus, an interesting problem of research has been to obtain
upper and lower bounds for $g(m)$.

In 1971, Erd\"{o}s and Selfridge [1, Theorem 3] proved
$g(m)\geq(1+o(1))\log m$. In 1975, as an improvement of results of
Cijsouw, Tijdeman [7] and Ramachandra [8], Ramachandra, Shorey and
Tijdeman [9] obtained an important result which states that $(\log
m/\log\log m)^3\ll g(m)$ by using Gelfond-Baker's theory. This
implies that Grimm's Conjecture would follow from Cram\'{e}r's
famous conjecture [3] which states that $p_{r+1}-p_r=O ((\log
p_r)^2)$. In 2006, Shanta and Shorey [10] confirmed that Grimm's
Conjecture is true for $m\leq 1.9\times 10^{10}$  and for all $n$ by
using Mathematica.

The object of this section is to study a stronger function $w(m)$
which is the largest integer $n$   such that the binomial
Coefficients  $\binom{m+n}{n}$ may be written as
$\binom{m+n}{n}=\prod_{i=1}^{i=n}a_i, a_i|(m+i),a_i\in
N,a_i>1,(a_i,a_j)=1,1\leq i\neq j\leq n$. By a result of the author
[11], we see that every binomial coefficient has the representation:
$\binom{m+n}{n}=\prod_{i=1}^{i=n}a_i, a_i|(m+i),a_i\in
N,(a_i,a_j)=1,1\leq i\neq j\leq n$. Naturally, we want to know what
condition $m$ satisfies so that each $a_i>1$ since it will imply
that  $\binom{m+n}{n}$ must have at least $n$ distinct prime factors
in this case. In this paper, we will prove the following theorem
without using Hall's theorem:

\vspace{3mm}\noindent {\bf Theorem 1:~~}%
When $m>\prod_{p\leq n}p^{[\log_pn]}$, the binomial coefficient
$\binom{m+n}{n}$ has the representation:
$$\binom{m+n}{n}=\prod_{i=1}^{i=n}a_i, a_i|(m+i),a_i\in
N,a_i>1,(a_i,a_j)=1,1\leq i\neq j\leq n.$$

Note that $w(m)\leq g(m)$, $n^{\pi (n)}>\prod_{p\leq
n}p^{[\log_pn]}$ and it is easy to show by [12] that $\pi (n)\log n
\geq n$ when $n\geq 17$. Hence we have $\log m\ll w(m)\ll
\sqrt{m/\log m}$. So, by this lower bound, we have obtained an
analogical result of the theorem 3 in [1].

In Section 3, we will give the proof of Theorem 1. In Section 4, we
will try to point out that it is not easy to improve the lower bound
of $w(m)$ to $(\log m)^2$.

\section{Proof of Theorem 1}
To prove the theorem, we need some lemmas which reflect that the
binomial coefficients have some fascinating and remarkable
arithmetic properties again.

\vspace{3mm}\noindent{\bf  Lemma 1:~~}%
Let $\binom{m+n}{n}$ be the binomial coefficient with $m,n\in N$.
Then for any prime $p$, $v_p(\binom{m+n}{n})\leq t$, where
$t=\max_{1\leq i\leq n}\{v_p(m+i)\}$.

\vspace{3mm}\noindent{\bf  Proof:~~}%
For the proof of Lemma 1, see [11].

\vspace{3mm}\noindent{\bf  Lemma 2:~~}%
If $m+i \not \in H_n $ for some $1\leq i\leq n$, then there is a
prime number $p$ such that $p^{v_p(m+i)}>n$. Moreover, $1\leq
v_p(\binom{m+n}{n})\leq v_p(m+i)$.

\vspace{3mm}\noindent{\bf  Proof:~~}%
By our assumption and the definition of $H_n$, clearly, there is a
prime $p$ such that $p^{v_p(m+i)}>n$. Note that $v_p(m+j)<v_p(m+i)$
when $j\neq i$. Otherwise, $p^{v_p(m+i)}|(j-i)$. But $|j-i|<n$ which
is impossible. Thus, $v_p(m+i)=\max_{1\leq j\leq n}\{v_p(m+j)\}$ and
$v_p(\binom{m+n}{n})\leq v_p(m+i)$ by Lemma 1. On the other hand, if
we write $m=p^{v_p(m+i)}x+y$, where $x, y\in N\cup \{0\}$ with
$0\leq y<p^{v_p(m+i)}$, then we have $p^{v_p(m+i)}|(y+i)$ and
$[\frac{m+n}{p^{v_p(m+i)}}]-[\frac{m}{p^{v_p(m+i)}}]=1$ since
$[\frac{n}{p^{v_p(m+i)}}]=0$ and $n+y\geq i+y$. Hence,
$\sum_{j=1}^\infty
([\frac{m+n}{p^j}]-[\frac{m}{p^j}]-[\frac{n}{p^j}])\geq 1$, and
$1\leq v_p(\binom{m+n}{n})$. This completes the proof of Lemma 2.

\vspace{3mm}  It is worthwhile pointing out that $H_n$ has many
interesting properties which should be given further consideration.
For example, $H_2$ is an empty set, $H_3=\{6\}$, $H_4=\{4,6,12\}$
and so on. $H_n\subseteq H_{n+1}$; $|H_n|=\prod_{p\leq n}(1+[\log_p
n])-1-\pi (n)$; $H_n\subseteq \psi (x,n)$, where $\psi (x,n)$ is the
set of n-smooth integers in $[1,x]$ with $x=\prod_{p\leq
n}p^{[\log_pn]}$. For some details on smooth integers which relate
to factorization of integers, see [13]. For applications of smooth
numbers to various problems in different areas of number theory, see
[15]-[18].

\vspace{3mm}\noindent{\bf  Lemma 3:~~}%
If $m+i \not \in H_n $ for every $1\leq i\leq n$, then
$\binom{m+n}{n}$ has the representation:
$$\binom{m+n}{n}=\prod_{i=1}^{i=n}a_i, a_i|(m+i),a_i\in
N,a_i>1,(a_i,a_j)=1,1\leq i\neq j\leq n.$$

\vspace{3mm}\noindent{\bf  Proof:~~}%
Write $B=\binom{m+n}{n}=\prod_{i=1}^{i=k}p_i^{e_i}$ by the
fundamental theorem of arithmetic, where $k$ is the number of
distinct prime factors of $B$. Clearly, for every prime factor $p_i$
of $B$, there must be a number $i_j$ with $1\leq i_j\leq n$, such
that $v_{p_i}(m+i_j)=\max_{1\leq r\leq n}\{v_{p_i}(m+r)\}$. By Lemma
1, we have $e_i\leq v_{p_i}(m+i_j)$ for every $1\leq i\leq n$.
Therefore, we can choose the number $m+i_j$ such that
$$\frac{(m+1)...(m+n)}{n!}=\frac{\frac{m+1}{p_i^{v_{p_i}(m+1)}}...\frac{m+i_j-1}{p_i^{v_{p_i}(m+i_j-1)}}\frac{m+i_j}{p_i^{v_{p_i}(m+i_j)-e_i}}\frac{m+i_j+1}{p_i^{v_{p_i}(m+i_j+1)}}...\frac{m+n}{p_i^{v_{p_i}(m+n)}}}{\frac{n!}{p_i^{v_{p_i}(n!)}}}.$$
Notice that $\frac{m+1}{p_i^{v_{p_i}(m+1)}}$, ...,
$\frac{m+i_j-1}{p_i^{v_{p_i}(m+i_j-1)}}$,
$\frac{m+i_j}{p_i^{v_{p_i}(m+i_j)-e_i}}$,
$\frac{m+i_j+1}{p_i^{v_{p_i}(m+i_j+1)}}$, ...,
$\frac{m+n}{p_i^{v_{p_i}(m+n)}}$, $\frac{n!}{p_i^{v_{p_i}(n!)}}$ are
all integers. Thus, by reduction, $B$ has the representation:
$$B=\prod_{i=1}^{i=n}a_i, a_i|(m+i),a_i\in N,(a_i,a_j)=1,1\leq
i\neq j\leq n.$$ On the other hand, since $m+i \not \in H_n $ for
every $1\leq i\leq n$, hence there is a prime number $p_i$ such that
$p_i^{v_{p_i}(m+i)}>n$. By Lemma 2, we have $v_{p_i}
(m+i)=\max_{1\leq r\leq n}v_{p_i} (m+r)$ and $p_i|a_i$. Note that
prime numbers $p_1,...,p_n$ are distinct. This completes the proof
of Lemma 3.

\vspace{3mm}\noindent{\bf  Proof of Theorem 1:~~} For every $1\leq
i\leq n$, we write
$m+i=(\prod_{p^{v_p(m+i)}>n}p^{v_p(m+i)})(\prod_{1\leq
p^{v_p(m+i)}\leq n}p^{v_p(m+i)})$. If $m+i\in H_n$ for some $1\leq
i\leq n$, then, we have $m+i=\prod_{1\leq p^{v_p(m+i)}\leq
n}p^{v_p(m+i)}\leq \prod_{1\leq p^{v_p(m+i)}\leq n}p^{[\log_pn]}$.
It is a contradiction by our assumption. Therefore, $m+i\not\in H_n$
for all  $1\leq i\leq n$, and by Lemma 3, we get the assertion of
Theorem 1.

\vspace{3mm}\noindent{\bf  Corollary 1:~~} If $2\leq n\leq 7$ and
$m+1,...,m+n$ are consecutive composite numbers, then
$\binom{m+n}{n}$ can be written as
$$\binom{m+n}{n}=\prod_{i=1}^{i=n}a_i, a_i|(m+i),a_i\in
N,a_i>1,(a_i,a_j)=1,1\leq i\neq j\leq n.$$ Namely, the conjecture in
[11] is true when $2\leq n\leq 7$.

\vspace{3mm}\noindent{\bf  Proof:~~} Clearly, Corollary 1 holds when
$2\leq n\leq 6$. So, we only consider the case of $n=7$. By Theorem
1, when $m>2^2\times 3\times 5\times 7=420$,  Corollary 1 holds.
When $m\leq 420$, by the table of prime numbers, the set of 7
consecutive composite numbers must be a subset of one of the
following 13 sets:

$\{90, 91, 92, 93, 94, 95, 96\}$,

$\{114, 115, 116, 117, 118, 119, 120, 121, 122, 123, 124, 125,
126\}$,

$\{182,183, 184, 185, 186, 187, 188, 189, 190\}$,

$\{200, 201, 202, 203, 204, 205, 206, 207, 208, 209, 210\}$,

$\{212, 213, 214, 215, 216, 217, 218, 219, 220, 221, 222\}$,

$\{242, 243, 244, 245, 246, 247, 248, 249, 250\}$,

$\{284, 285, 286, 287, 288, 289, 290, 291, 292\}$,

$\{318, 319, 320, 321, 322, 323, 324, 325, 326, 327, 328, 329,
330\}$,

$\{338, 339,340, 341, 342, 343, 344, 345, 346\}$,

$\{360, 361, 362, 363, 364, 365, 366\}$,

$\{390, 391, 392, 393, 394. 395, 396\}$,

$\{402, 403, 404, 405, 406, 407, 408\}$,

$\{410, 411, 412, 413, 414, 415, 416, 417, 418\}$.

We have $\frac{204\times 205\times 206\times 207\times 208\times
209\times 210}{1\times 2\times 3\times 4\times 5\times 6\times 7
}=17\times 41 \times 103\times 207\times 208 \times 209\times 5$.
The remaining cases satisfy $m+i\not\in H_7$ for every $1\leq i\leq
7$ when $n=7$. By Lemma 3, Corollary 1 holds.

\section{Remarks}
In this section, we will try to explain that it is not easy to
improve the lower bound of  $w(m)$ to $(\log m)^2$. On one hand, we
found two exceptions to the conjecture in [11] when $m=116,n=10$ or
$m=118,n=8$. This implies that the conjecture in [11] is not always
true. On the other hand, if the lower bound of  $w(m)$  can be
improved to $(\log m)^2$, then there are only finitely many
exceptions to the conjecture in [11] which implies Grimm's
conjecture (Note that Grimm Conjecture follows from Cram\'{e}r's
conjecture by using results of Ra-machandra, Shorey and Tijdeman.).
Namely, one would guess the following:

\vspace{3mm}\noindent{\bf  Conjecture 1:~~} For every sufficiently
large integer $n$, the product of $n$  consecutive composite numbers
$m+1,...,m+n$ may be written
$$\prod_{i=1}^{i=n}(m+i)=n!\times \prod_{i=1}^{i=n}a_i, a_i|(m+i),a_i\in
N,a_i>1,(a_i,a_j)=1,1\leq i\neq j\leq n.$$

Obviously, if Conjecture 1 holds, then it also implies a surprising
consequence that there are only finitely many exceptions to the
following Conjecture 2.

\vspace{3mm}\noindent{\bf  Conjecture 2:~~} (i) For every $n>1$,
there are primes in the intervals $(d-n,d+n]$ and $(d-q,d+q)$
respectively, where $d>1$ is a factor of $\prod_{n<p<2n}p$,  and $q$
is the least prime factor of $d$.

(ii) For every $n>2$, let $d>1$ be a factor of $\prod_{n/2<p<n}p$.
If $d=nt+r$  and $d$ is coprime to each of
$nt+1,...,nt+r-1,nt+r+1,...,nt+n$, then there are primes in the
interval $[nt+1,nt+n]$.

\vspace{3mm}\noindent{\bf  Proof:~~} (i) By Bertrand-Chebyshev
theorem which states that there exists a prime in interval $(n, 2n)$
when $n > 1$,we see that $\prod_{n<p<2n}p$ has a prime factor $>n$.
Let $d>1$ be a factor of $\prod_{n<p<2n}p$. If there are no primes
in the interval $(d-n,d+n]$, then $d-n+1,...,d-n+2n$ are all
composite numbers. By the assumption that Conjecture 1 holds, we
have that $\prod_{i=1}^{i=2n}(d-n+i)=(2n)!\times
\prod_{i=1}^{i=2n}a_i$ with $a_i|(d-n+i),a_i\in
N,a_i>1,(a_i,a_j)=1,1\leq i\neq j\leq 2n $ for  sufficiently large
integer $n$. But this is impossible since we obtain $a_n=1$ from
$v_p((2n)!)=1$ and $v_p((d-n+1)...(d-n+2n))=1$.

  Let $q$ be the least prime factor of $d$. If $d$ is prime, then
  $q=d$. Clearly, there are primes, for example $q$, in the interval
$(d-q,d+q)=(0,2q)$. Now, we assume that $d$ is not prime. If there
are no primes in the interval $(d-q,d+q)$, then $d-q+1,...,d-q+2q-1$
are all composite numbers. Conjecture 1 holds, so we have that
$\prod_{i=1}^{i=2q-1}(d-q+i)=(2q-1)!\times \prod_{i=1}^{i=2q-1}a_i$
with $a_i|(d-q+i),a_i\in N,a_i>1,(a_i,a_j)=1,1\leq i\neq j\leq 2q-1
$ for sufficiently large integer $n$. Note that for any prime
divisor $r$ of $\frac{d}{q}$, $r>q$ since $d>1$ is a factor of
$\prod_{n<p<2n}p$ and is square-free. Notice also that $r<2q$ since
$q>n$ and $r<2n$. Therefore $a_q=1$ since $v_r((2q-1)!)=1$ and
$v_r((d-q+1)...(d-q+2q-1))=1$. This contradiction shows that the
case (i) holds.

Similarly, one could deduce the case (ii) in Conjecture 2 assuming
Conjecture 1. In fact, if there are not any primes in the interval
$[nt+1,nt+n]$, then by the assumption that Conjecture 1 holds, we
have that $\prod_{i=1}^{i=n}(nt+i)=n!\times \prod_{i=1}^{i=n}a_i$
with $a_i|(nt+i),a_i\in N,a_i>1,(a_i,a_j)=1,1\leq i\neq j\leq n $
for every sufficiently large integer $n$. However,  by our
assumption and the known conditions, for any prime divisor $p$ of
$d=nt+r$, we have $v_p(n!)=1$ and $v_p((nt+1)...(nt+n))=1$ since
$d=nt+r$ is square-free and  coprime to each of
$nt+1,...,nt+r-1,nt+r+1,...,nt+n$. Therefore $a_r=1$ since
$v_p(\binom{nt+n}{n})=0$ which is a contradiction.

\vspace{3mm} Conjecture 2 implies also that for every $n>1$, there
are primes in the interval $(d-q,d+q)$, where $d=\prod_{n<p<2n}p$,
and $q$ is the least prime factor of $d$. Of course, if $d$ is
prime, then there are primes in the interval $(d-q,d+q)=(0,2q)$.
  When $d$ is not prime,  or equivalently, when $n=4$ or $n>5$ (by the refined Bertrand-Chebyshev
theorem which states that there exists at least two distinct primes
in interval $(n, 2n)$ when $n=4$ or $n>5$), we have $d-q\neq 0$.
Moreover, by Chebyshev theorem which states that $\frac{0.92129
x}{\log x}<\pi (x)<\frac{1.1056x}{\log x}$, we have that
$2n<2q<4n<((\frac{0.92\times 2n}{\log 2n}-\frac{1.11n}{\log n})\log
n)^2<((\pi (2n)-\pi (n))\log q)^2<(\log(\prod_{n<p<2n}p-q))^2$ for
every sufficiently large integer $n$. Thus, $2q<(\log(d-q))^2$ which
is out of bounds for even the Cram\'{e}r's conjecture, and it is
hard to prove that there are primes in the interval $(d-q,d+q)$, say
nothing of Conjecture 1. Therefore, we think that it is not easy to
improve the lower bound of $w(m)$ to $(\log m)^2$. Is Conjecture 1
true?

\vspace{3mm}Based on Conjecture 2, one could obtain some interesting
results. For instance, the below is a fast algorithm for generating
large primes by using small known primes.

\vspace{3mm}\noindent{\bf An Algorithm for generating an m-bit prime
number~~}

\vspace{3mm}\textbf{Input:} A natural number $m$

\textbf{Step1:} Choose $p_1<...<p_r<2p_1$ such that $r$ is
appropriately large and $p_i$ are all prime (need not consecutive),
for $1\leq i\leq r$ and $2^{m-1}<k=\prod_{i=1}^{i=r}p_i<2^m$.

This step can be finished easily by pre-computing.

\textbf{Step2:} Test whether each of $k\pm 2,...,k\pm 2[p_1/2]$ is
prime. If for some $i$, $k+i$ or $k-i$ is prime, terminate the
algorithm. Otherwise, Conjecture 2 does not hold for $n=2[p_1/2]+1$.

Such a prime can be found quickly when $r$ is appropriately large.
On one hand, primality testing is comparatively easy since its
running time is polynomial [19]. On the other hand, since $k$ has
many prime divisors, hence there is a high probability that either
$k+2i$ or $k-2i$ is prime for some $1\leq i\leq
[p_1/2]<p_1/2<2^{\frac{m}{r}-1}$.

\textbf{Output:} An m-bit prime number

\vspace{3mm} Due to the fact that it lies outside the scope of this
paper. We omitted more details of this algorithm.  As a toy example,
we can generate a 32-bit prime using small primes 29, 31, 37, 41,
43, 47. Note that $2^{31}<29\times 31\times 37\times 41\times
43\times 47<2^{32}$. By an exhaustive search, or a simple sieve, we
find that $29\times 31\times 37\times 41\times 43\times
47+6=2756205449$ is a 32-bit prime. Unfortunately, we do not know
whether Conjecture 2 is true or not, although we can test whether
Conjecture 2 holds by the aforementioned algorithm. Moreover, we
have not been able to work out a complete proof of Conjecture 2,
still less conjecture 1. But, all these and related questions,
specially, the lower bound of  $w(m)$ we hope to investigate.

\section{Acknowledgements}
I am happy to record my gratitude to my advisor Professor Xiaoyun
Wang who made a lot of helpful suggestions. I wish to thank the
referees for a careful reading of this paper. Thank the Institute
for Advanced Study in Tsinghua University for providing me with
excellent conditions. This work was partially supported by the
National Basic Research Program (973) of China (No. 2007CB807902)
and the Natural Science Foundation of Shandong Province (No.
Y2008G23).

\clearpage
\end{document}